\documentclass[runningheads]{llncs}

\usepackage{url}
\makeatletter
\AtBeginDocument{%
  \def\doi#1{\upshape\url{https://doi.org/#1}}}
\makeatother

\usepackage{graphicx}
\usepackage{lipsum}
\usepackage{amsfonts}
\usepackage{graphicx}
\usepackage{epstopdf}
\usepackage{algorithmic}
\usepackage{tabularx}
\usepackage{mathtools}
\usepackage{amsopn}
\DeclareMathOperator{\diag}{diag}
\DeclareMathOperator{\sinc}{sinc}
\DeclareMathOperator*{\argmin}{argmin}

\DeclarePairedDelimiter{\norm}{\|}{\|}
\DeclarePairedDelimiter{\abs}{\lvert}{\rvert}
\DeclarePairedDelimiter{\lidx}{[}{]}

\graphicspath{{./}{./graphics}}

\begin{document}
\title{Adaptive Regularization of B-Spline Models for Scientific Data%
    \thanks{This work is supported by the U.S. Department of Energy, Office of Science, Advanced Scientific Computing Research under Contract DE-AC02-06CH11357, Program Manager Margaret Lentz.}%
}
%
%
\author{David Lenz\inst{1}
    \and Raine Yeh\inst{2}
    \and Vijay Mahadevan\inst{1}
    \and Iulian Grindeanu\inst{1}
    \and Tom Peterka\inst{1}
}
\authorrunning{D. Lenz et al.}
%
\institute{Argonne National Laboratory, Lemont, IL, USA \\
        \email{\{dlenz, mahadevan, iulian, tpeterka\}@anl.gov}
    \and Google Inc., New York, NY, USA \\
        \email{raineyeh@google.com}
}

\maketitle              
\begin{abstract}
    B-spline models are a powerful way to represent scientific data sets with a functional approximation. However, these models can suffer from spurious oscillations when the data to be approximated are not uniformly distributed. Model regularization (i.e., smoothing) has traditionally been used to minimize these oscillations; unfortunately, it is sometimes impossible to sufficiently remove unwanted artifacts without smoothing away key features of the data set. In this article, we present a method of model regularization that preserves significant features of a data set while minimizing artificial oscillations. Our method varies the strength of a smoothing parameter throughout the domain automatically, removing artifacts in poorly-constrained regions while leaving other regions unchanged. The behavior of our method is validated on a collection of two- and three-dimensional data sets produced by scientific simulations.

\end{abstract}

\section{Introduction}
Data sets assembled from scientific simulations or experimental readings are often defined as a list of position-value pairs, where each data point consists of a measurement and the corresponding location of that measurement. These point locations can form structured grids, unstructured meshes, or unconnected point clouds, depending on the application. Methods for analyzing these data often apply only to particular layouts; usually, numerical analysis techniques become more complex as the geometry of the point locations becomes more general (e.g. point clouds). Even seemingly straightforward tasks such as interpolation can be computationally burdensome on unstructured point clouds and numerically inaccurate on highly nonuniform meshes. One way to avoid these challenges is by representing a data set with a mathematical function and then analyzing the function instead of the original data.  This can substantially streamline the process of interpolation and differentiation away from data points, simplify visualization tasks, and make resampling the data almost trivial.

The focus of this article is the approximation of scientific data sets by (tensor-product) B-splines. B-splines are a family of smooth functions used ubiquitously throughout geometric modeling~\cite{lin2018geometric} and form the underpinnings of isogeometric analysis (IGA)~\cite{hughes2005iga}. Recent study has shown that large, complex data sets produced by scientific simulations at extreme scale can be effectively modeled by B-splines~\cite{peterka_ldav18}. Similar results have also been obtained for nonuniform rational B-splines, or NURBS, which are a generalization of B-splines~\cite{nashed2019}.

B-splines have a number of properties that make them useful as a functional representation of data. 
B-splines are high-order approximants, and evaluating, differentiating, and integrating a B-spline model is fast and numerically stable~\cite{deBoor1978guide}. Crucially, differentiation and integration can be computed in closed-form and incur no additional loss of accuracy, unlike finite differences or Riemann sums. Thus, once a sufficiently accurate spline has been computed to represent the data, it is often more productive to analyze the functional model than the original data.

However, computing a best-fit B-spline requires solving a linear system which may be ill-conditioned or rank-deficient. A common cause of this ill-conditioning is an input data set that contains both sparse and dense patches of points in proximity to each other. 
Without additional effort, solving this system can produce a function that oscillates strongly between data points or even diverges in regions where input data are very sparse. This problem can be hard to detect automatically, since error metrics are usually defined in terms of the pointwise error between the original data and the model. Spurious oscillations occurring away from the input data will not be captured by these metrics.

To address these challenges, we developed a new method for fitting B-splines to unstructured data that reduces or eliminates oscillations while leaving critical features of the data set unchanged. Our method regularizes the solution to the B-spline fitting problem by adding a variable-strength smoothing parameter that automatically adapts based on characteristics of the input data set. This additional term smooths out spike artifacts in regions where the data set is very sparse but does not do any smoothing where data points are densely packed, thereby preserving accuracy in these regions. In addition, our method creates well-defined spline models even for data sets with irregular boundaries. No knowledge of the boundary is required; the method automatically handles areas outside the boundary that contain no data points. 

The remainder of this paper is organized as follows. A review of related ideas and methods is given in Section~\ref{sec:related-work}. In Section~\ref{sec:background}, we provide a primer on the mathematical details used to describe B-splines throughout the paper. Our main result, a method for adaptive regularization of B-spline models, is described in Section~\ref{sec:method}. We then exhibit the performance of this method in Section~\ref{sec:results} with a series of numerical examples. We summarize directions for further research in Section~\ref{sec:future-work} and present conclusions in Section~\ref{sec:conclusions}.

\section{Related Work}\label{sec:related-work}
Creating B-spline models to represent unstructured data sets is a particular example of scattered data approximation (SDA), a broad area of study concerned with defining continuous functions that interpolate or approximate spatially scattered inputs. SDA is often applied to image reconstruction problems, where an experimental or physical constraint prohibits the collection of uniformly-spaced samples, such as medical~\cite{arigovindan2006full}, seismic~\cite{duijndam1999reconstruction}, or astronomical~\cite{vio2000reconstruction} imaging. An introductory comparison of SDA methods was compiled by Francis et al.~\cite{francis2018scattered}.

Ill-conditioned numerical methods are a persistent challenge throughout the SDA literature, and a number of techniques have been proposed to increase numerical stability. Our approach is most similar to the variational methods of SDA, in which the magnitude of the approximating function's derivative (or ``roughness'') is minimized. The early work of Duchon~\cite{duchon1977splines} is a canonical example. Historically, roughness minimization has been achieved through the use of smoothing splines; a thorough exposition of smoothing splines can be found in the book by Gu~\cite{gu2013smoothing}. The application of smoothing splines requires a trade-off between accuracy and roughness minimization, since aggressively penalizing roughness tends to degrade accuracy. Therefore, much work has been devoted to parametrizing this trade-off appropriately. Craven and Wahba~\cite{craven1978smoothing} developed the influential ``cross-validation'' approach, which  is expanded upon by Gu~\cite{gu1992cross}.

The functional approximation used in this article is based on global tensor-product B-splines, but a number of other spline-based regression methods have been proposed. Truncated thin-plate splines were used by Wood~\cite{wood2003thin} to improve the efficiency of thin-plate regression splines while maintaining their characteristic stability. Lee et al.~\cite{lee1997scattered} utilized hierarchies of B-splines to fit unstructured data points, but the instability arising from sparse point distributions was not treated explicitly. Francis et al.~\cite{francis2018scattered} consider a two-step process for resampling unstructured point clouds with variable point density onto unstructured grids. While this method does not construct a functional approximation, it does show good performance as a resampling methodology.

Our novel adaptive regularization procedure was first explored in the dissertation of the second author~\cite{yehDissertation}. The method is also directly inspired by the work of El-Rushaidat et al.~\cite{elrushaidat2021}, in which a two-level regularization process was introduced in the context of resampling unstructured data onto structured meshes. However, their method requires an ad-hoc selection of the criteria to switch between high and low regularization strengths, as well as an application-dependent overall level of smoothing. A notable contribution in our work is a continuously varying regularization strength (not two-level) that is adapted automatically.

\section{Background on B-Splines}\label{sec:background}
In this section, we provide a brief overview of the basic definitions and constructions necessary to describe B-spline models for scientific data. A thorough presentation on the fundamental theory of B-splines can be found in the books by de Boor~\cite{deBoor1978guide} and by Piegl and Tiller~\cite{piegl1997nurbs}.
\subsection{B-Spline Curves}
A one-dimensional B-spline curve of degree $p$ in $\mathbb{R}^D$ is a parameterized curve
\begin{equation}\label{eq:spline-def}
C(u) = \sum_{j=0}^{n-1} N_{j,p}(u) P_j,
\end{equation}
where each $N_{j,p}$ is a piecewise-polynomial function of degree $p$, and each $P_j \in \mathbb{R}^D$ is a ``control point'' in $D$-dimensional space. 

The B-spline basis functions, denoted $N_{j,p}$, are defined on the parameter space $[0,1] \subset \mathbb{R}$, which is divided by a nondecreasing sequence of ``knots'' $t_0 \leq t_1 \leq \ldots \leq t_{n+p} \in [0,1]$. Each basis function $N_{j,p}$ is a bump function in $[t_j, t_{j+p+1}]$ and zero elsewhere.\footnote{%
  We consider only ``clamped'' knot sequences in this paper; thus, the first $p+1$ knots are always $0$ and the last $p+1$ knots are always $1$.%
}
In this paper, we will assume that the degree of the B-spline is fixed and drop the $p$ subscript, instead denoting the $j^{th}$ function as $N_j$.

In order to simplify notation when describing high-dimensional tensor product splines, we use multi-indices to index quantities in multiple dimensions simultaneously.  A multi-index $\alpha = (\alpha^1, \ldots, \alpha^d)$ is a $d$-tuple of nonnegative indices, where the sum of components of $\alpha$ is denoted $\abs{\alpha} = \sum_k \alpha_k$. 

We will often consider index sets for our multi-indices in the form of 
\begin{equation}\label{eq:index-set}
  A = \{\text{all } \alpha \in \mathbb{N}^d \text{ such that } 0 \leq \alpha^k < n_k \text{ for } 1 \leq k \leq d\},
\end{equation}
where $n_k$ are previously defined positive numbers. We impose a lexicographic ordering on these sets, 
and denote by $\lidx{\alpha}_A$ the index of $\alpha$ in the lexicographic ordering of $A$.
In the following sections, we consider matrices in which each column corresponds to a multi-index. In this scenario, we list multi-indices in lexicographic order; thus, the multi-index $\alpha$ corresponding to the $j^{th}$ column satisfies $\lidx{\alpha}_A = j$.

\subsection{Tensor Product B-Splines}
Tensor product B-splines are a natural extension of B-spline curves to higher-dimensional manifolds, such as surfaces, volumes, and hypervolumes. Here, we denote by $d$ the dimension of the tensor product volume and $D$ the dimension of the ambient space (for instance, a 2D surface in 3D space would correspond to $d=2$, $D=3$). The parameter space for a $d$-dimensional tensor product B-spline is $[0,1]^d$, which is divided by $d$ different knot vectors $\mathbf{t}_k = \{t_k^j\}_{j=0}^{n_k + p}$, $k=1,\ldots,d$. 

Given a tuple $u = (u^1, \ldots, u^d) \in [0,1]^d$, the tensor product basis functions are defined as
\begin{equation}
  N_\alpha(u) = \prod_{k=1}^d N^k_{\alpha^k}(u^k).
\end{equation}
where $\alpha$ is a multi-index as described above and $N_{\alpha^k}^k$ is the $(\alpha^k)^{th}$ basis function with respect to the knot vector $\mathbf{t}_k$.
With $n_k + p + 1$ total knots in each dimension, there are $n_k$ basis functions in each dimension.\footnote{%
Here we assume for simplicity that the degree of the B-spline is the same in each dimension, but the degree can vary in practice if desired.}
Therefore, the total number of tensor product basis functions is $n_{tot} = \prod_{k=1}^d n_k$, which is also the total number of control points for the the tensor product spline.

A $d$-dimensional tensor product B-spline in $\mathbb{R}^D$ is a function of the form
\begin{equation}\label{eq:tensor-spline}
  C(u) = \sum_{\alpha \in A} N_\alpha(u) P_\alpha,
\end{equation}
where $A$ is the set of all basis functions and $P_\alpha \in \mathbb{R}^D$ for each $\alpha$.

\subsection{Optimal Control Points}
Given a collection of knot vectors and polynomial degree, the best-fit B-spline to a given data set is determined by a linear least-squares minimization problem. Let $\{Q_i\}_{i=0}^{m-1}$ be the list of points in $\mathbb{R}^D$ to be approximated with a $d$-dimensional tensor product spline. For each $0 \leq i < m$, let $v_i \in [0,1]^d$ be the parameter tuple corresponding to the point $Q_i$. The optimal control points are determined by the least-squares minimization problem:
\begin{equation}\label{eq:standard-min}
  \{\hat{P}_j\} = \argmin_{P_j} \sum_{i=0}^{m-1} \norm{Q_i - C(v_i)}^2.
\end{equation}

This minimization problem can be rewritten in normal form by differentiating the objective function in Equation~\eqref{eq:standard-min} with respect to each of the control points. The normal system reduces to the matrix equation $\mathbf{N}^T \mathbf{N} \mathbf{P} = \mathbf{N}^T \mathbf{Q}$, where
\begin{equation}
  \mathbf{N}_{ij} = N_\alpha(v_i) \quad \text{where } \lidx{\alpha}_A = j, \qquad
  \mathbf{P}_{ij} = P_\alpha^j, \quad \text{where } \lidx{\alpha}_A = i, \qquad
  \mathbf{Q}_{ij} = Q_i^j.
\end{equation}
The superscripts in the above equations index the components of the vectors $P_\alpha$ and $Q_i$.
$\mathbf{N}$ is a $m \times n_{tot}$ matrix, often called the ``B-spline collocation matrix,'' $\mathbf{P}$ is an $n_{tot} \times D$ matrix with each row containing a control point, and $\mathbf{Q}$ is a $m \times D$ matrix with each row containing an input point.

Typically, the matrix $\mathbf{N}$ is very sparse, and this system may be solved with an iterative method or sparse direct solver. However, as we show in the following section, this system is ill-conditioned when the sample density of the input points $P_i$ varies from region to region.

\section{Adaptive Regularization}\label{sec:method}
\begin{figure}
  \centering
  \includegraphics[width=110pt, height=60pt]{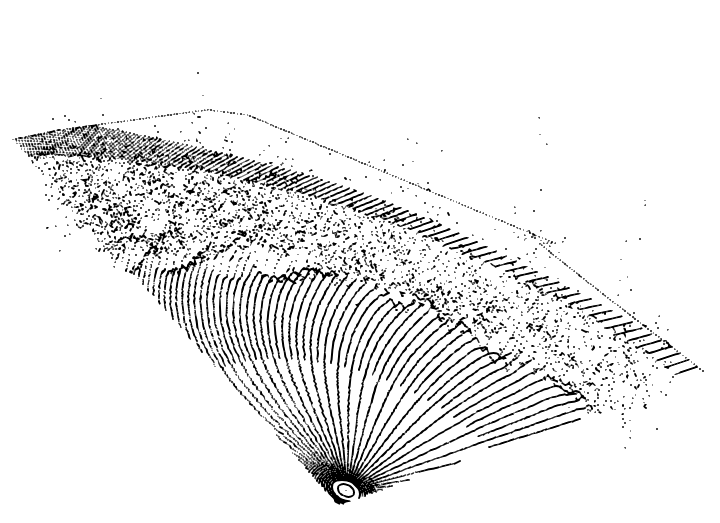}
  \includegraphics[width=110pt, height=60pt]{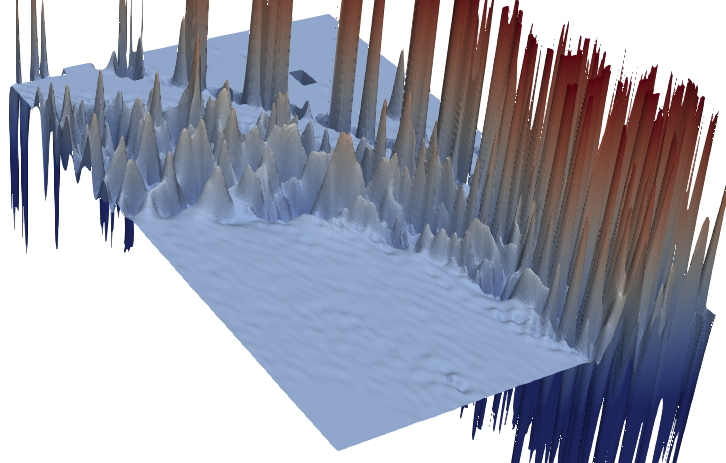}
  \includegraphics[width=110pt, height=60pt]{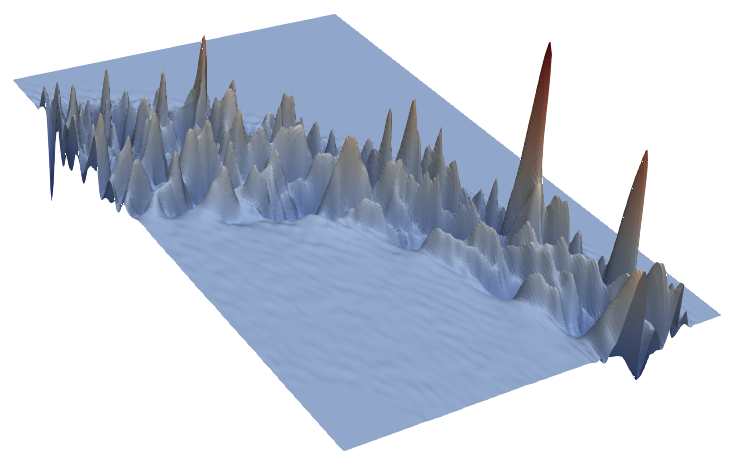}
  \caption{Left: A data set with nonuniform point density. Center: Best-fit B-spline model without regularization. Right: B-spline model with adaptive regularization. The center image is cropped; spike artifacts in this model extend well outside the frame.}
  \label{fig:xgc-sample}
\end{figure}

A significant challenge when modeling unstructured data with tensor product B-splines is the (ill-)conditioning of the fitting procedure. Generally speaking, tensor product B-spline models can oscillate strongly due to overfitting in regions where input data is sparse (see Figure~\ref{fig:xgc-sample}).
Our method of adaptive regularization produces a unique solution to systems which would otherwise be rank-deficient and improves the overall conditioning of the system. In practice, the adaptively regularized models possess fewer oscillatory artifacts and do not exhibit any divergent behavior in our testing. In contrast to standard regularization techniques, our method does not smooth out the model indiscriminately -- instead, it regularizes only those regions that require smoothing.

This technique employs a spatially-varying regularization strength that is computed automatically as a function of the relative positioning of input data points to the B-spline knots. In general, the regularization strength increases in regions with little to no input data and decreases (potentially to zero) in regions ``saturated'' with input points. When the regularization strength is zero throughout a region of the domain, no smoothing is performed in that region; therefore, any sharp features present in densely sampled regions of the domain will be preserved by the adaptive regularization procedure.

Standard roughness minimization can be formulated as a penalized least-squares minimization problem, similar to Equation~\eqref{eq:standard-min}. The penalty term weights the size of the second derivative at each point with a new parameter, which we denote as $\lambda > 0$. The control points of the regularized spline are defined by:
\begin{equation}\label{eq:penalized-min}
  \{\hat{P}_j\}  = \argmin_{P_j} \left(\sum_{i=0}^{m-1} \norm{Q_i - C(v_i)}^2 + \lambda^2 S(C) \right),
\end{equation}
where $S(C)$ approximates the size of the second derivatives of $C$.

Let $w_\alpha \in [0,1]^d$ be the parameter that maximizes the value of $N_\alpha$. Let $\partial^\delta$ denote the partial derivative where the order of derivative in each dimension is given by the components of multi-index $\delta$.\footnote{%
For example, $\partial^{(2,0)} f = \partial^2 f / \partial x_1^2$, and $\partial^{(0,2)} f = \partial^2 f / \partial x_2^2$, while $\partial^{(1,1)} f = \partial^2 f / (\partial x_1 \partial x_2)$.
} 
We define $S(C)$ to be 
\begin{equation}\label{eq:deriv-approx}
  S(C) = \sum_{\alpha \in A} \sum_{\abs{\delta}=2} \norm*{\partial^{\delta}C(w_\alpha)}^2.
\end{equation}
Note that the summation above is a sum over all derivatives of order 2, including mixed partial derivatives.

Equation~\eqref{eq:penalized-min} can be converted into a system of equations in the same way as Equation~\eqref{eq:standard-min}. The only additional step is computing the derivative of $S(C)$ with respect to the control points. 
In matrix form, the system is 
\begin{equation}
  \begin{pmatrix} \mathbf{N}^T & \lambda \mathbf{M}^T \end{pmatrix}
  \begin{pmatrix} \mathbf{N} \\ \lambda \mathbf{M} \end{pmatrix} \mathbf{P} = 
  \mathbf{N}^T \mathbf{Q}, \qquad \text{where }
  \mathbf{M} = 
  \begin{pmatrix} \mathbf{M}_{\delta_1} \\ \vdots \\ \mathbf{M}_{\delta_n} \end{pmatrix}
\end{equation}
and $(\mathbf{M}_{\delta})_{i,j} = \partial^\delta N_\beta(w_\alpha)$, $\lidx{\alpha}_A = i$, $\lidx{\beta}_A = j$. Intuitively, each column of $\mathbf{M}_\delta$ describes the $\partial^\delta$ partial derivative of an individual B-spline basis function. The matrix $\mathbf{M}$ is the concatenation of all the individual $\mathbf{M_\delta}$ matrices, where $\abs{\delta} = 2$. $\mathbf{N}$, $\mathbf{P}$, and $\mathbf{Q}$ are defined as in Section~\ref{sec:background}. 

The novel improvement of our adaptive regularization scheme is to modify the above system of equations by varying the size of $\lambda$ for each column of $\mathbf{M}$. Since each column of this matrix corresponds to a B-spline basis function and control point, variation in the size of $\lambda$ provides a mechanism to modify the smoothing conditions on each control point of the spline individually.
Due to the local support property of B-splines, setting $\lambda_j = 0$ for control points in a given region ``disables'' the regularization in that region, while still allowing for smoothing to be applied elsewhere.
Algebraically, we replace the scalar parameter $\lambda$ by a diagonal matrix $\mathbf{\Lambda} = \diag(\lambda_1,\ldots, \lambda_{n_{tot}})$, where each $\lambda_j \geq 0$, and consider the new linear system
\begin{equation}\label{eq:penalized-matrix}
  \begin{pmatrix} \mathbf{N}^T & (\mathbf{M}\mathbf{\Lambda})^T \end{pmatrix}
  \begin{pmatrix} \mathbf{N} \\ \mathbf{M}\mathbf{\Lambda} \end{pmatrix} \mathbf{P} = 
  \mathbf{N}^T \mathbf{Q}.
\end{equation}
The value of each $\lambda_i$ is computed automatically as a function of the relative positioning between input data points and B-spline knots. 

To better control this function, we introduce a user-specified parameter called the ``regularization threshold,'' denoted $s^*$. Changing the regularization threshold adjusts the criterion by which some regions of the domain are smoothed and others are not. As $s^*$ increases, smoothing constraints will be applied to larger and larger regions in the domain.

Let $s_j$ denote the $j^{th}$ column sum of $\mathbf{N}$ and $\widetilde{s}_j$ the $j^{th}$ column sum of $\mathbf{M}$. Given $s^* \geq 0$, we define
\begin{equation}\label{eq:lambda-j}
  \lambda_j = \frac{\max(s^* - s_j,0)}{\widetilde{s}_j}.
\end{equation}
Thus, $\mathbf{\Lambda}$ is defined such that every column sum of $\binom{\mathbf{N}}{\mathbf{M\Lambda}}$ is no less than $s^*$.

Adapting the regularization strengths $\lambda_j$ this way has a number of important results. 
If $s^* = 0$, then $\mathbf{\Lambda} = 0$ and the minimization becomes the usual least-squares problem. When $s^*$ is small, $\lambda_j$ will be zero unless the $j^{th}$ column sum of $\mathbf{N}$ is small, which is indicative of an ill-conditioned system. Here, adaptive regularization smooths out only those control points which are poorly constrained. 

This formulation also explains why the adaptive regularization method preserves sharp, densely sampled features while smoothing out oscillatory artifacts. In regions of the domain that are densely sampled, control points will be constrained by many data points and thus the corresponding column sum in $\mathbf{N}$ will be relatively large. By choosing $s^*$ to be sufficiently small, all control points in this region will have a regularization strength of zero; i.e. $\lambda_j = 0$. Therefore, the best-fit spline in this region will not be artificially smoothed.

Finally, we remark that the adaptive regularization framework, while defined above in terms of second derivatives, can easily be extended to other derivatives. We find that constraining first and second derivatives simultaneously is particularly helpful when modeling data sets with no points at all in certain regions.  This typically happens when the data represent an object with an interior hole or irregular boundary. To consider both first and second derivatives, the only change is to matrix $\mathbf{M}$ in Equation~\eqref{eq:penalized-matrix}. Originally, $\mathbf{M}$ is the concatenation of all matrices $\mathbf{M}_\delta$, where $\delta$ describes a second derivative. To minimize first derivatives as well, we change this to the concatenation of all $\mathbf{M}_\delta$ such that $\delta$ describes a first or second derivative. Both versions of the method are described in Section~\ref{sec:results}.

\section{Results}\label{sec:results}
We demonstrate the effectiveness of our method with a series of numerical experiments. First, we compare adaptive regularization to uniform regularization where the smoothing parameter has been chosen manually. Next, we study the reconstruction of an analytical signal from sparse samples with varying levels of sparsity. 
For each sparsity level, we report the error and condition number for and unregularized and adaptively regularized model.
 We then test the performance of adaptive regularization on data sets with no data in certain regions. In these problems, we construct a B-spline model that extrapolates into regions with no pointwise constraints, and check that the adaptive regularization method produces a reasonable result.

\subsection{Data Sets}
The performance of the adaptive regularization method was studied on a collection of two- and three-dimensional point clouds with different characteristics. Some data sets were sampled from analytical functions so that we could compute pointwise errors relative to a ground truth, while other data sets were generated by scientific experiments and simulations.

\textit{2D Polysinc.}
The polysinc data set is a two-dimensional point cloud sampling the function
$f(x,y) = \sinc\left(x^2 + y^2\right) \sinc\left(2(x-2)^2 + (y+2)^2\right)$.\footnote{%
  We consider the unnormalized sinc function: $\sinc(x) = \sin(x)/x$, with $\sinc(0)=1$.}
360,000 point locations are uniformly sampled from the box domain $[-4\pi,4\pi] \times [-4\pi,4\pi]$, except at four disk-shaped regions where the sample rate is 50$\times$ lower. 

\textit{XGC Fusion.}
The XGC fusion data set represents a normalized derivative of electrostatic potential in a single poloidal plane of a Tokamak fusion simulation. The data set contains 56,980 points with an irregular boundary and was produced by the XGC code~\cite{xgc} in a simulation of the gyrokinetic equations.

\textit{CMIP6 Climate.}
The CMIP6 climate data set represents ocean surface temperature in a projected box region around Antarctica. The data set contains 585,765 points with a large hole (representing Antarctica) in the center and was produced by a Coupled Model Intercomparison Project (CMIP6)~\cite{cmip6} simulation.

\textit{sahex Nuclear.}
The sahex nuclear data set is derived from a simulation of a single nuclear reactor component, produced with the SHARP toolkit~\cite{yu2016sharp}. The three-dimensional data are bounded by a hexagonal prism and point density is coarser in the $z$ dimension than $x$ and $y$. The data set contains 63,048 points.

\subsection{Comparison of Adaptive versus Uniform Regularization}
\begin{figure}
  \centering
  \vspace{-20pt}
  \begin{minipage}{.42\textwidth}
    \includegraphics[width=\textwidth]{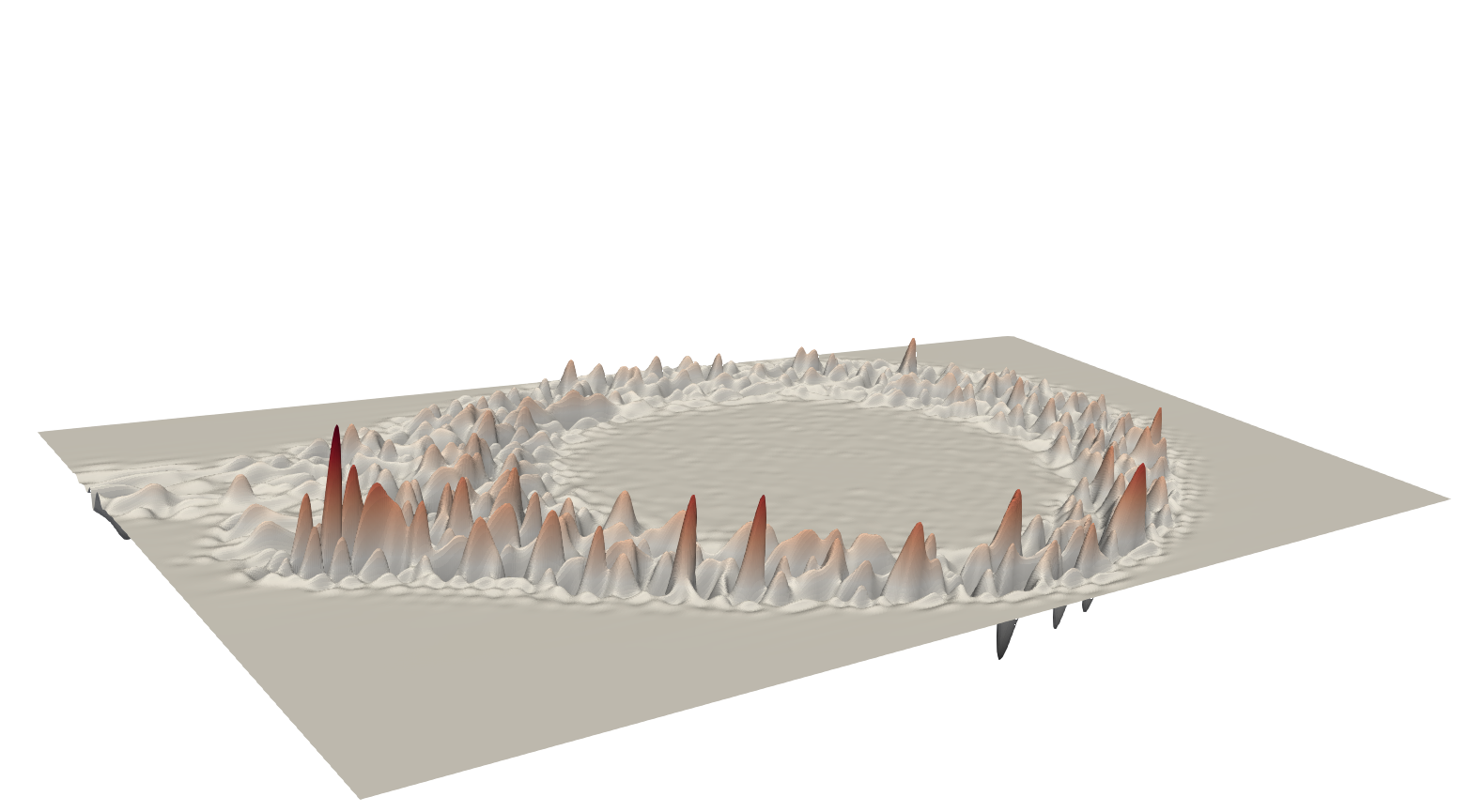}\\[-35pt]
    \includegraphics[width=\textwidth]{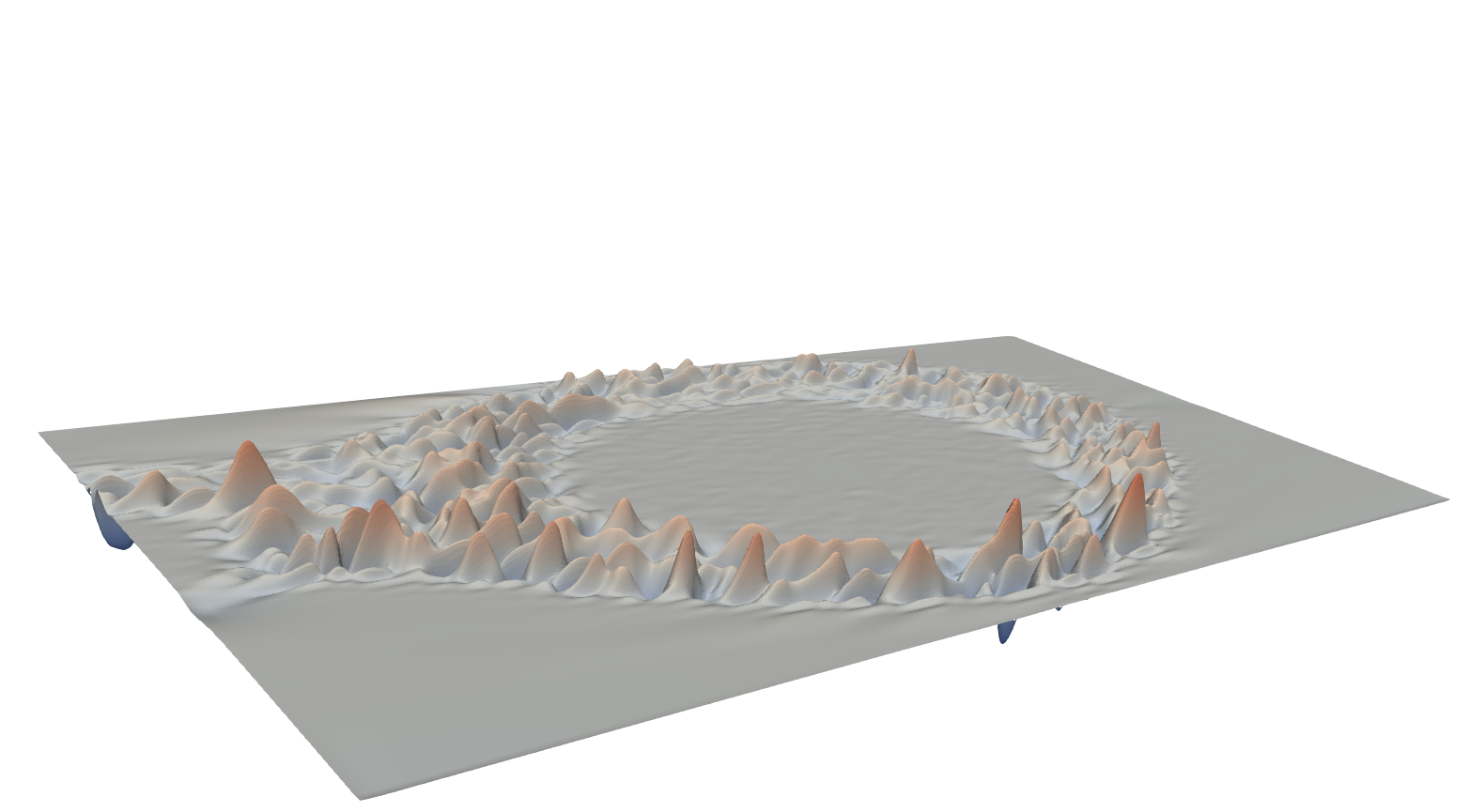}
  \end{minipage}
  \begin{minipage}{.42\textwidth}  
    \includegraphics[width=\textwidth]{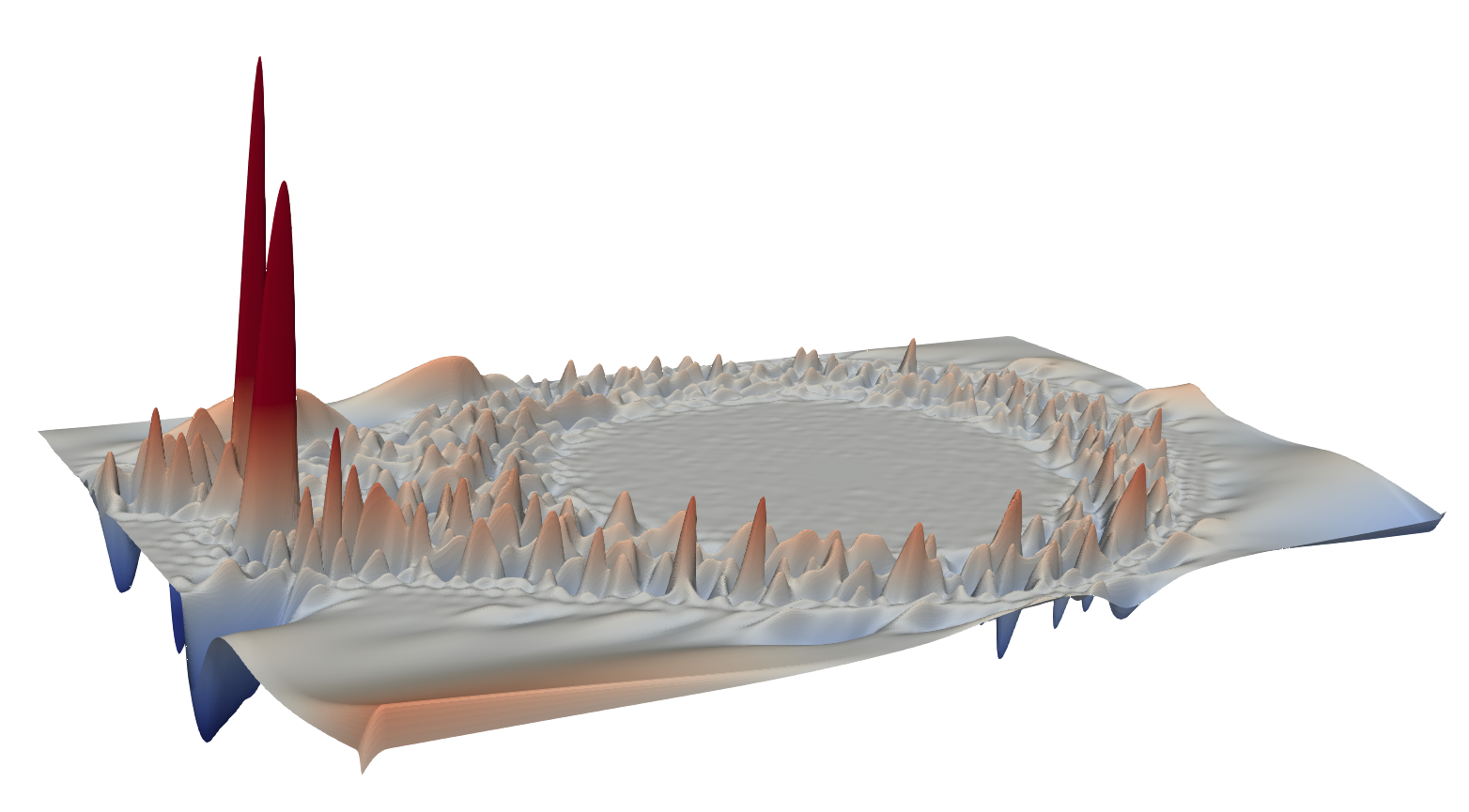}\\[-35pt]
    \includegraphics[width=\textwidth]{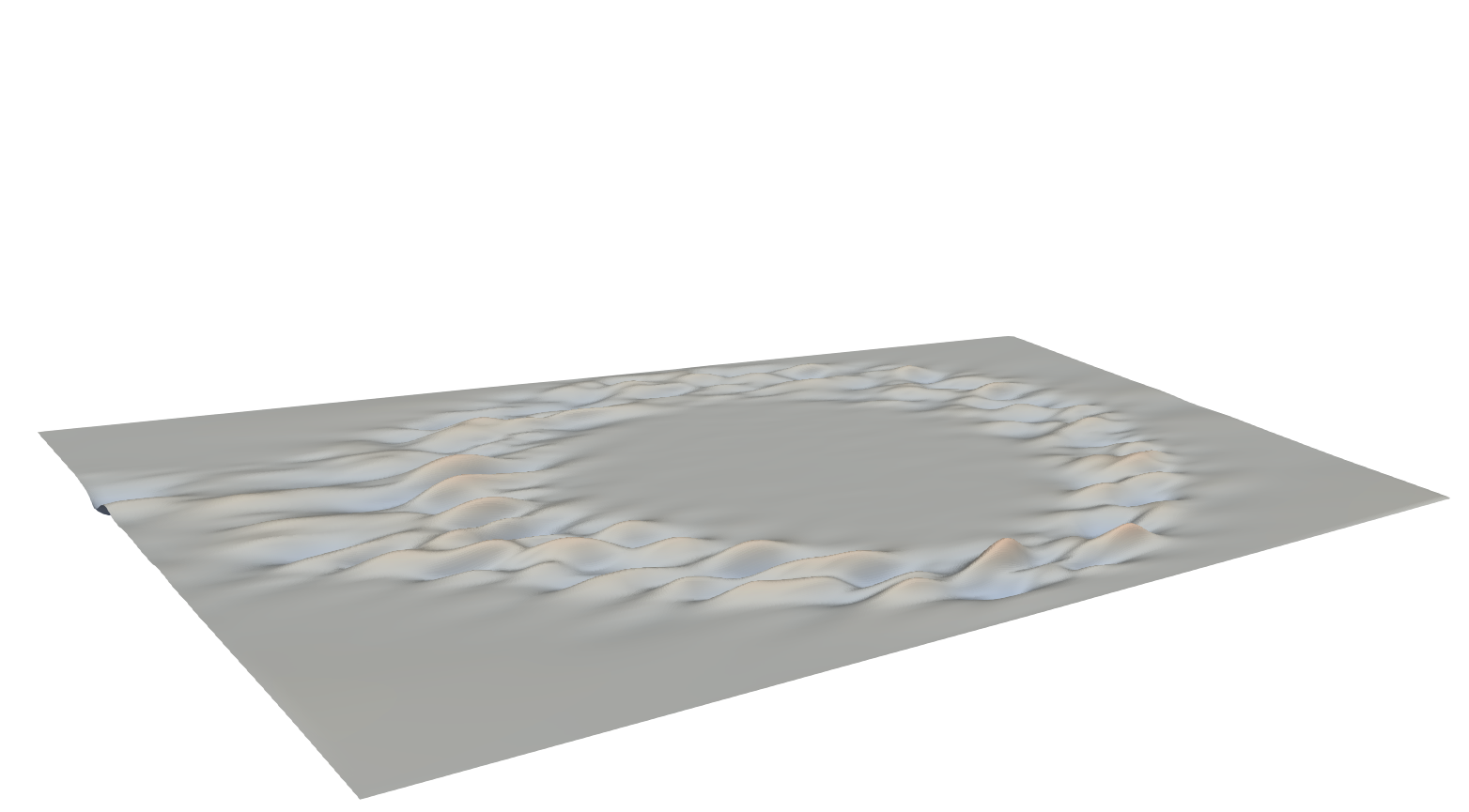}
  \end{minipage}
  \raisebox{-60pt}{\includegraphics[width=30pt]{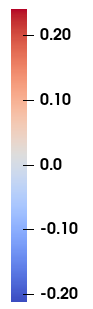}}
  \caption{Comparison of uniform vs adaptive regularization. Clockwise from top-left: Adaptive regularization, uniform regularization with $\lambda = 10^{-6}$, uniform regularization with $\lambda = 10^{-4}$, uniform regularization with $\lambda = 10^{-5}$.}
  \label{fig:uniform-vs-adaptive} 
\end{figure}

Applying a uniform regularization strength to an entire model can produce unsatisfactory results, because sufficiently smoothing oscillatory artifacts can also smooth out sharp features. Figure~\ref{fig:uniform-vs-adaptive} compares our adaptive regularization scheme (with $s^*=6$) against three strengths of uniform regularization. 
The data in Figure~\ref{fig:uniform-vs-adaptive} is the XGC fusion data set, which contains sharp peaks in a ring but is flat inside the ring. Data are sparse or nonexistent outside the ring. The two images at right show a model with uniform regularization that is too weak, causing artifacts (top), or too strong, dampening the features (bottom). The best uniform regularization strength we could find is given at bottom-left, but even in this example the characteristic peaks in the data are smoothed down.

\vspace{-4ex}

\subsection{Accuracy on Analytical Signals}
\begin{figure}
  \centering
  \includegraphics[height=80pt, width=110pt]{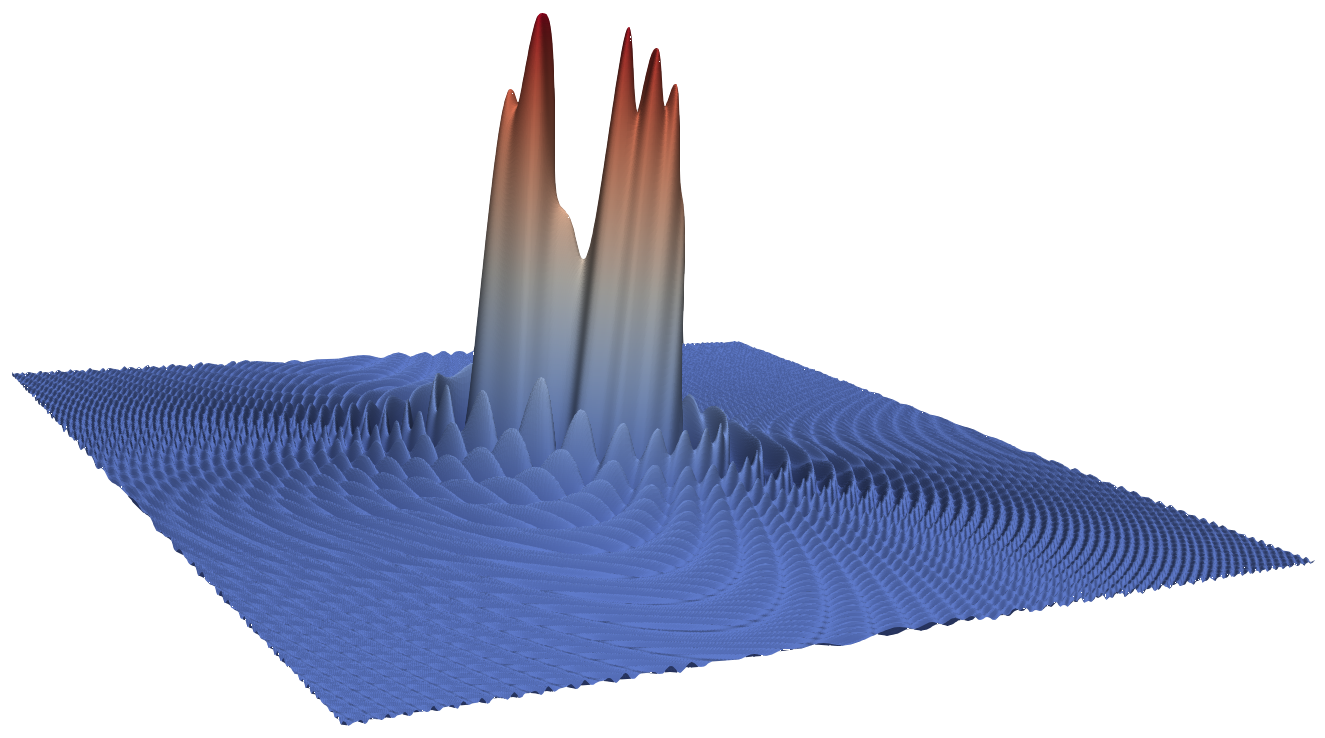} 
  \includegraphics[height=80pt, width=110pt]{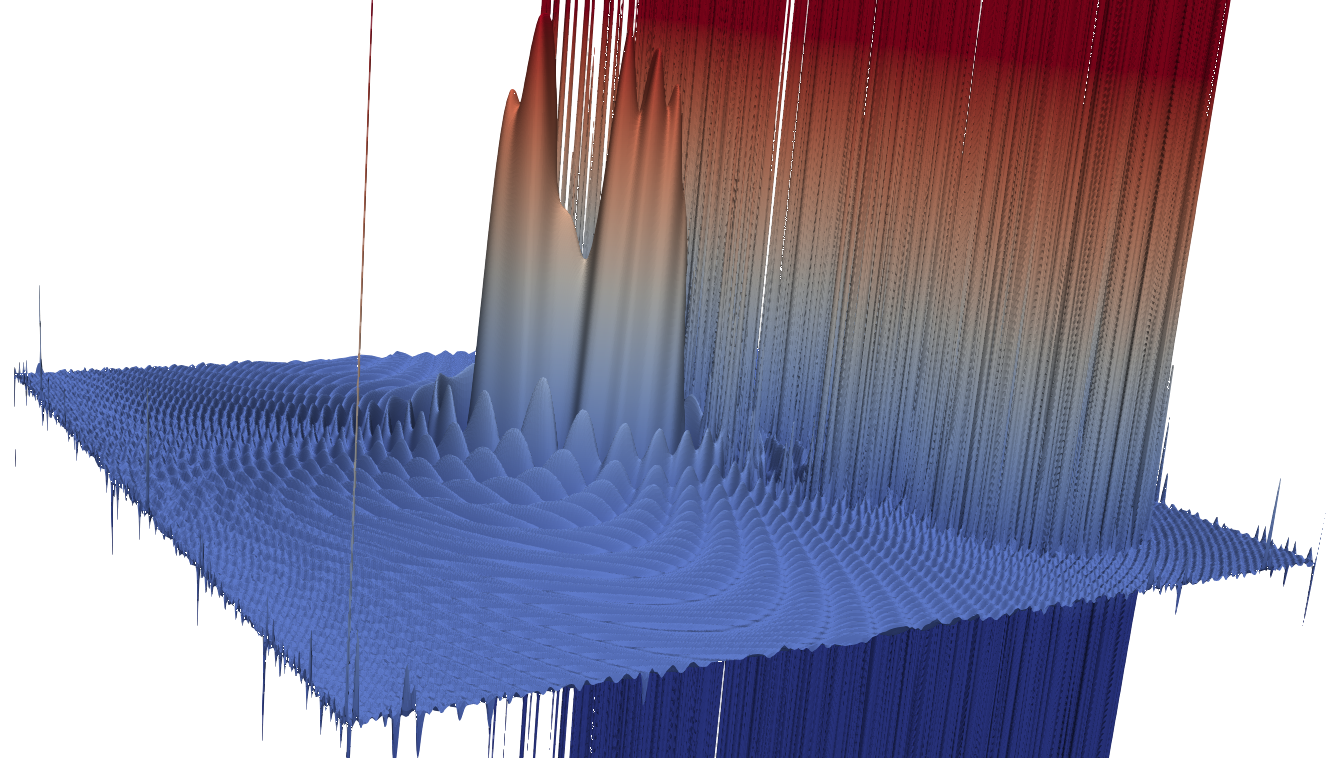}
  \includegraphics[height=80pt, width=110pt]{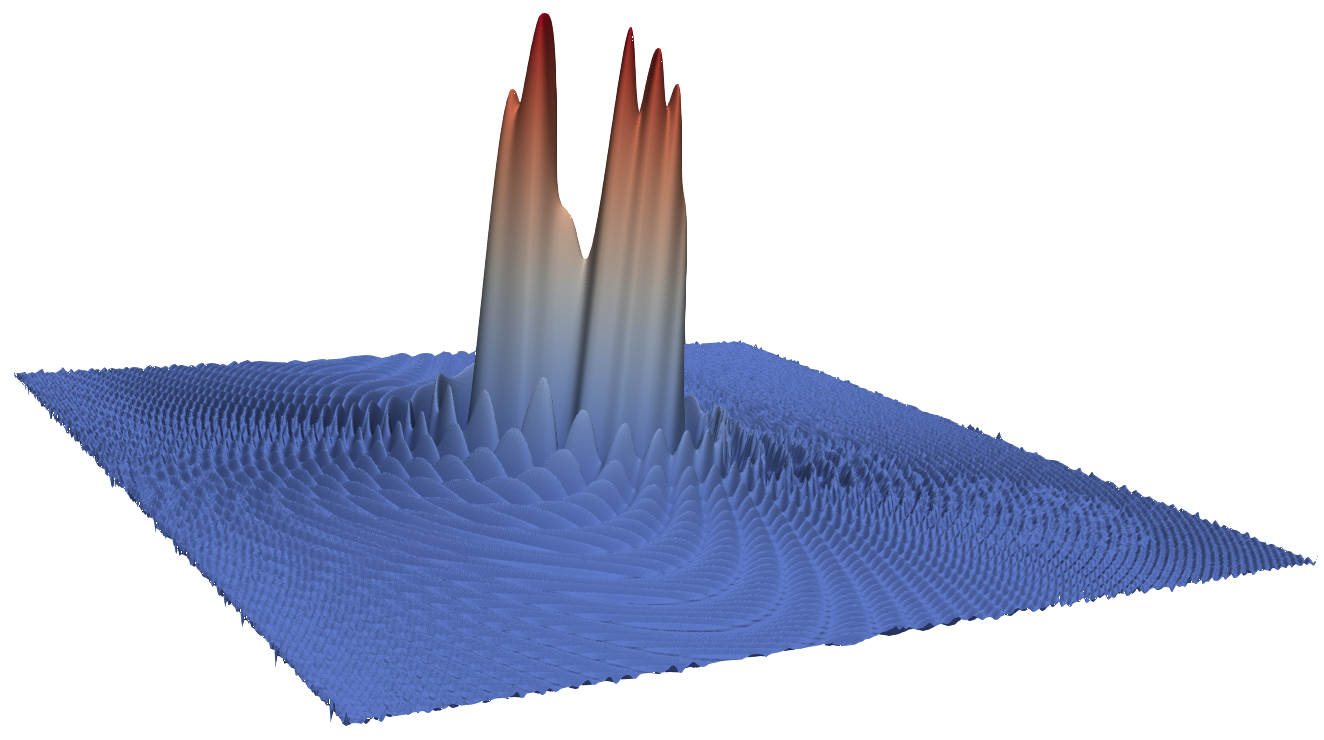}\\
  \hspace{25pt}
  \includegraphics[height=80pt, width=90pt]{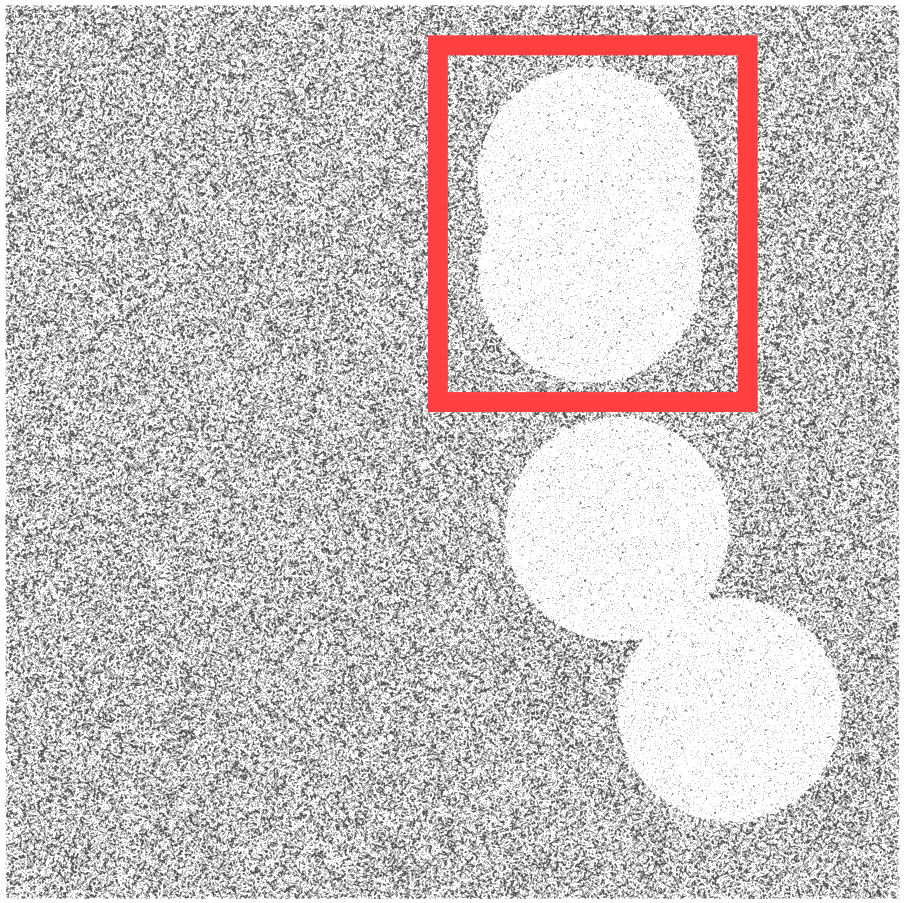}
  \hspace{5pt}
  \includegraphics[height=80pt, width=90pt]{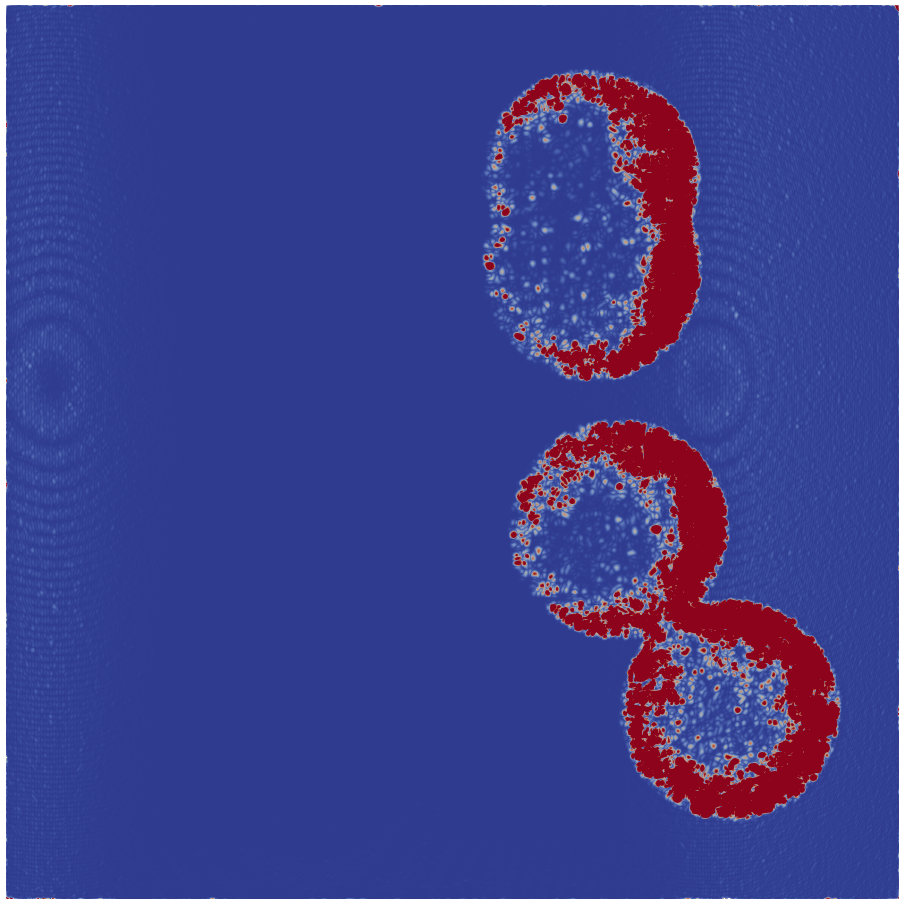}
  \hspace{5pt}
  \includegraphics[height=80pt, width=90pt]{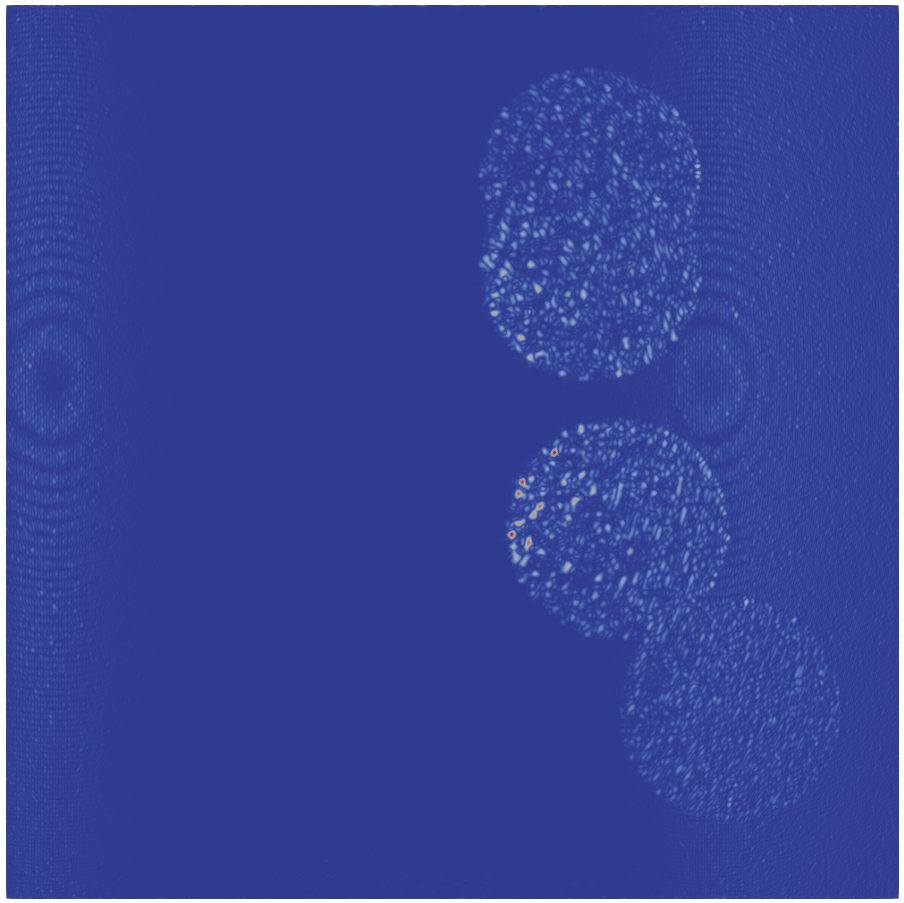}
  \includegraphics[width=25pt]{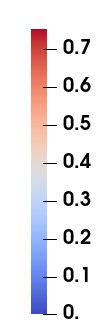}
  \caption{Top row: Synthetic polysinc signal (left), model with no regularization (center), model with adaptive regularization (right). Bottom row: Top down view of input distribution (left), error profile with no regularization (center), error profile with adaptive regularization (right). Area of interest for error calculation is in red at bottom-left.}
  \label{fig:psinc-2d}
\end{figure}

To quantify the accuracy of B-spline models with adaptive regularization, we consider the oscillatory polysinc function with a highly nonuniform input data set (Figure~\ref{fig:psinc-2d}). We illustrate two B-spline models, one fit without regularization and one with our adaptive regularization ($s^*=1$). Both models are degree four with a 300 $\times$ 300 grid of control points. Without regularization, the model diverges in the regions of low sample density; with adaptive regularization, the model produces an accurate representation even where sample density is low. A top-down view of the error profiles is given in the second row of Figure~\ref{fig:psinc-2d}. A close comparison of the ground-truth (top left) and adaptively-regularized spline (top right) shows that the spline model is not artificially smoothed in dense regions, even though the sparse regions are smoothed.  In particular, our regularization procedure preserves the distinctive oscillations and peaks in the signal.

\vspace{-4ex}
The degree of sparsity in the input data strongly influences the accuracy of a B-spline model. Table~\ref{tb:error-study} lists the errors in each model for varying levels of sparsity in the voids. The errors are measured in a box around two voids (see Figure~\ref{fig:psinc-2d}) in order to pinpoint the behavior of the models in this region. When the point density is equal inside and outside of the voids (sparsity$=1.0$), error for both models is low. As the voids become more sparse, the error in the unregularized model increases by four orders of magnitude while error in the adaptively regularized model stays essentially flat.

\begin{table}[h!]
  \caption{Model errors as a function of sparseness. Maximum and $L^2$ (average) errors are computed for both adaptively regularized and unregularized models in the vicinity of two voids (see Figure~\ref{fig:psinc-2d}). Condition numbers for both minimization problems are reported at bottom.}
\centering
\begin{tabularx}{\textwidth}{l @{\hspace{5pt}}|@{\hspace{3pt}} 
  *6{>{\centering\arraybackslash}X}}
Sparsity & 0.02 & 0.08 & 0.16 & 0.32 & 0.64 & 1.00\\ \hline
Max Error (reg) & 3.25e-2 & 2.89e-2 & 2.71e-2 & 2.39e-2 & 1.20e-2 & 1.16e-2 \\
Max Error (no reg)  & 1.29e2 & 6.27e2 & 2.36e-1 & 3.65e-2 & 1.20e-2 & 1.16e-2 \\ \hline
$L^2$ Error (reg) & 1.93e-3 & 1.53e-3 & 1.13e-3 & 7.04e-4 & 5.56e-4 & 5.44e-4 \\
$L^2$ Error (no reg) & 2.17e0 & 4.68e0 & 3.42e-3 & 7.35e-4 & 5.56e-4 & 5.44e-4\\ \hline
Condition \# (reg) & 177 & 980 & 289 & 198 & 121 & 189 \\
Condition \# (noreg) & inf & inf & 6.97e4 & 2.58e3 & 1.57e3 & 3.74e3 
\end{tabularx}
\label{tb:error-study}
\end{table}

Data sparsity also affects the condition number of the least-squares minimization. Table~\ref{tb:error-study} gives the condition number of the matrices $\mathbf{N}$ (`noreg') and $\binom{\mathbf{N}}{\mathbf{M\Lambda}}$ (`reg') for each sparsity level. As sparsity is increased, the condition number of $\binom{\mathbf{N}}{\mathbf{M\Lambda}}$ remains lower and steady but the condition number of $\mathbf{N}$ starts higher and eventually becomes infinite. Condition numbers were computed with the Matlab routine \texttt{svds}.

\subsection{Extrapolation into Unconstrained Regions}
When a large region of the domain does not contain any data points to constrain the best-fit B-spline problem, the least-squares minimization will be ill-posed and the resulting model can exhibit extreme oscillations. However, data sets with empty regions or ``holes'' are very common in scientific and industrial applications. For example, some climate models measure ocean temperatures or land temperatures, but not both simultaneously. Industrial simulations often model objects with irregular boundaries, and data from physics simulations are shaped by the locations of detectors.

Although empty regions are usually omitted in subsequent analysis, it is still important to understand and control the behavior of a B-spline model in empty regions. Extreme oscillations near the boundary of a hole can distort the derivative of the model away from this boundary. In addition, attempting to compute simple statistics about the model (minimum, maximum, mean) can be biased if the model exhibits unpredictable behavior in empty regions.

\begin{figure}
  \centering
  \hspace{36pt}
  \begin{minipage}{0.25\textwidth}
    \includegraphics[width=\textwidth]{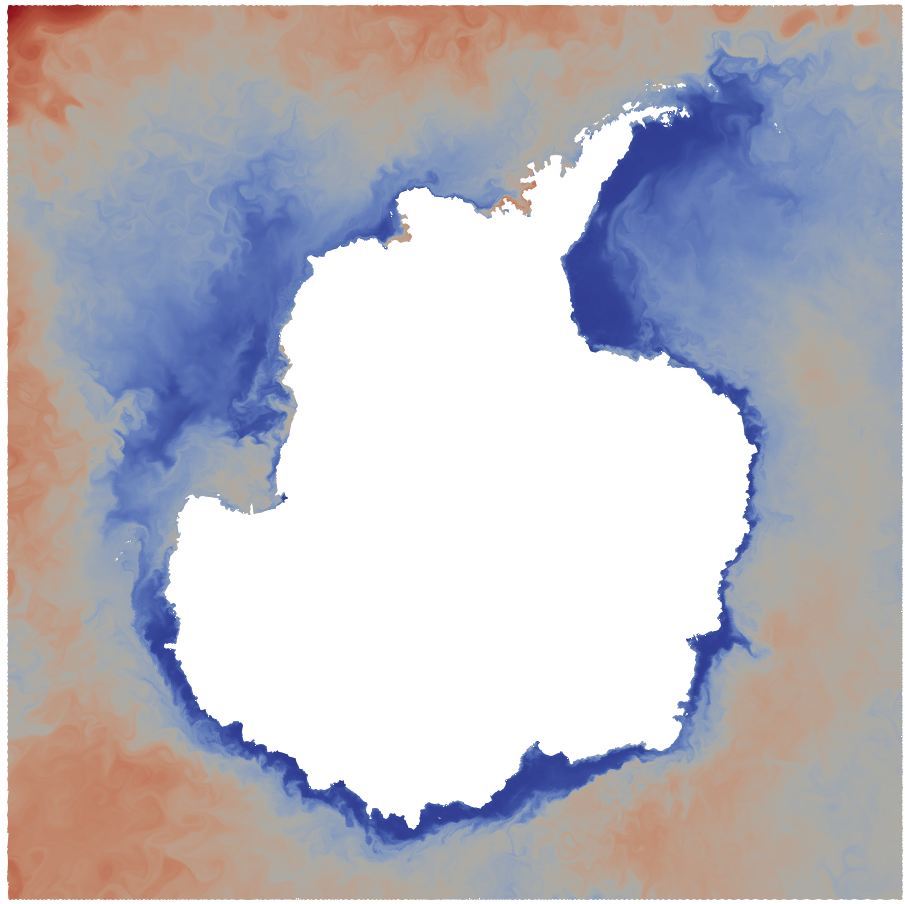} 
    \includegraphics[width=\textwidth]{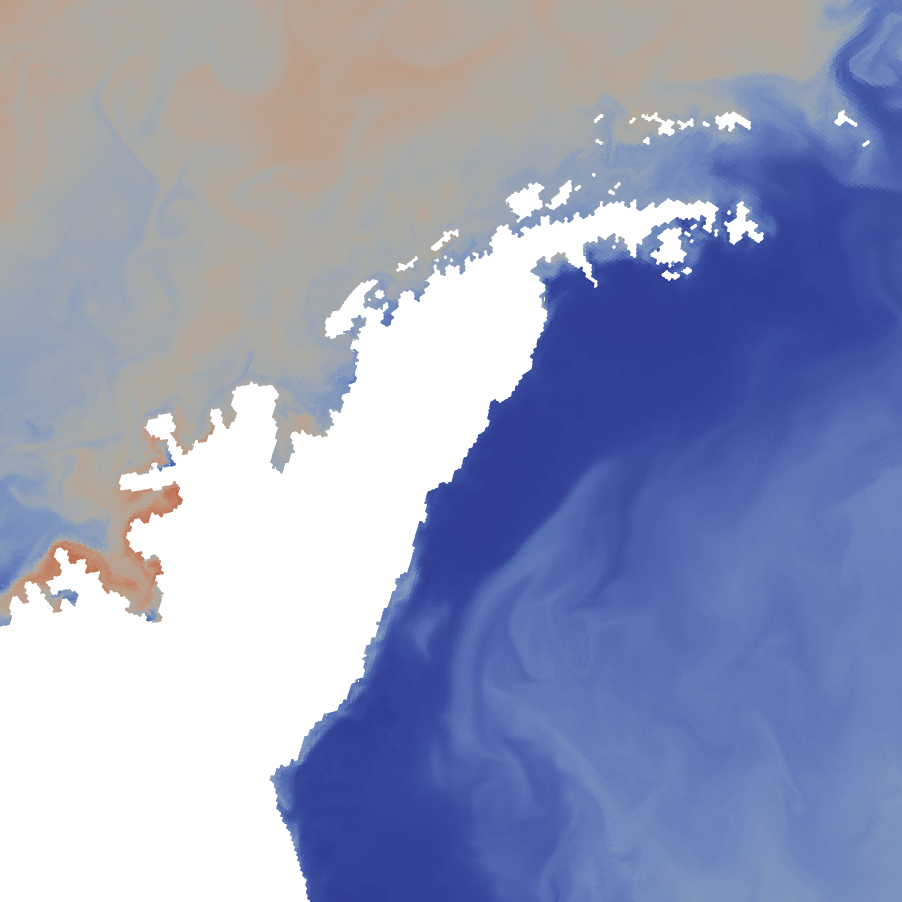} 
  \end{minipage}
  \begin{minipage}{0.25\textwidth}
    \includegraphics[width=\textwidth]{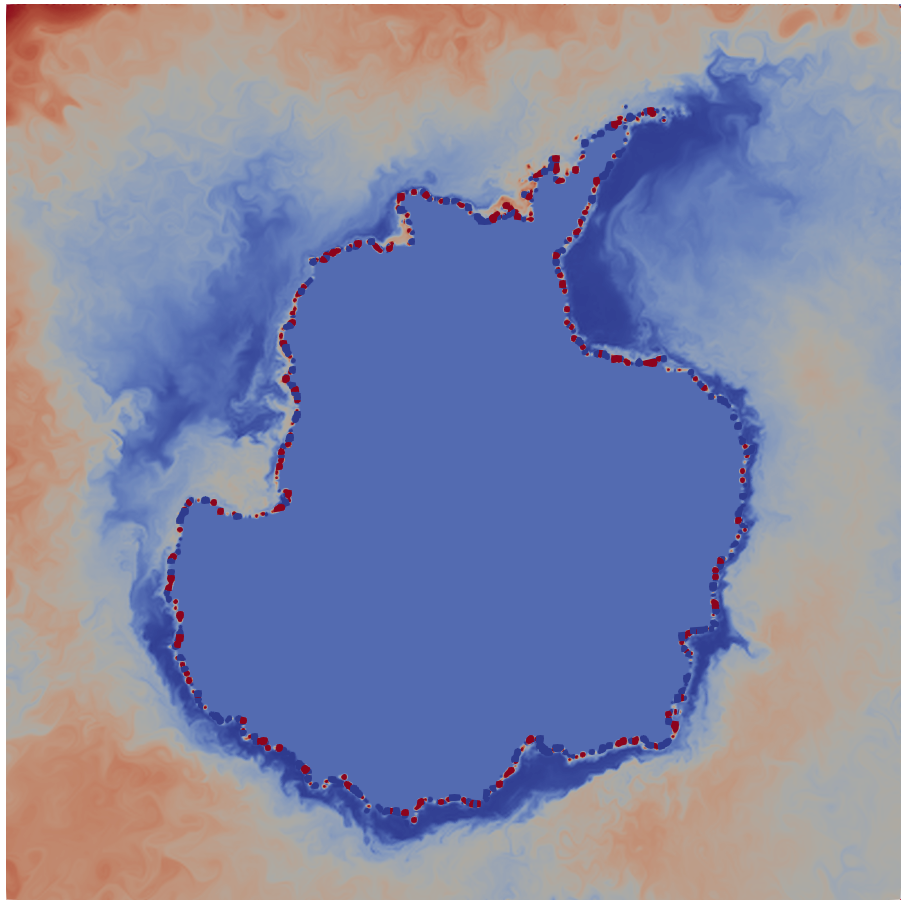}
    \includegraphics[width=\textwidth]{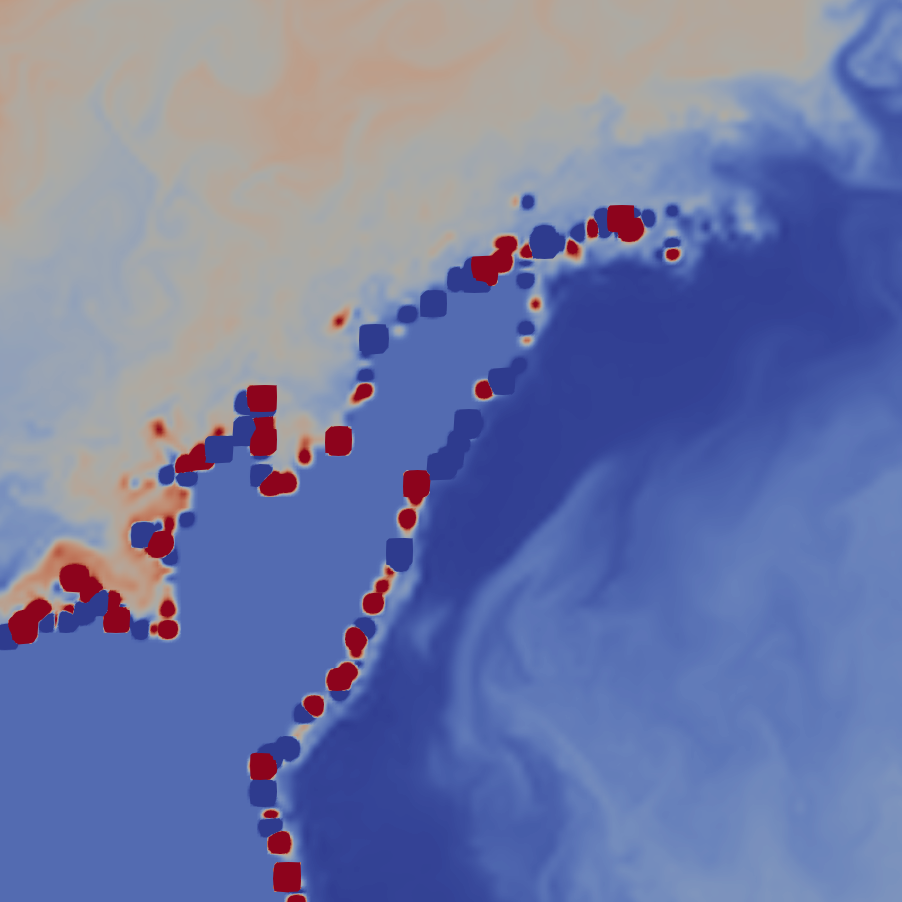}
  \end{minipage}
  \begin{minipage}{0.25\textwidth}
    \includegraphics[width=\textwidth]{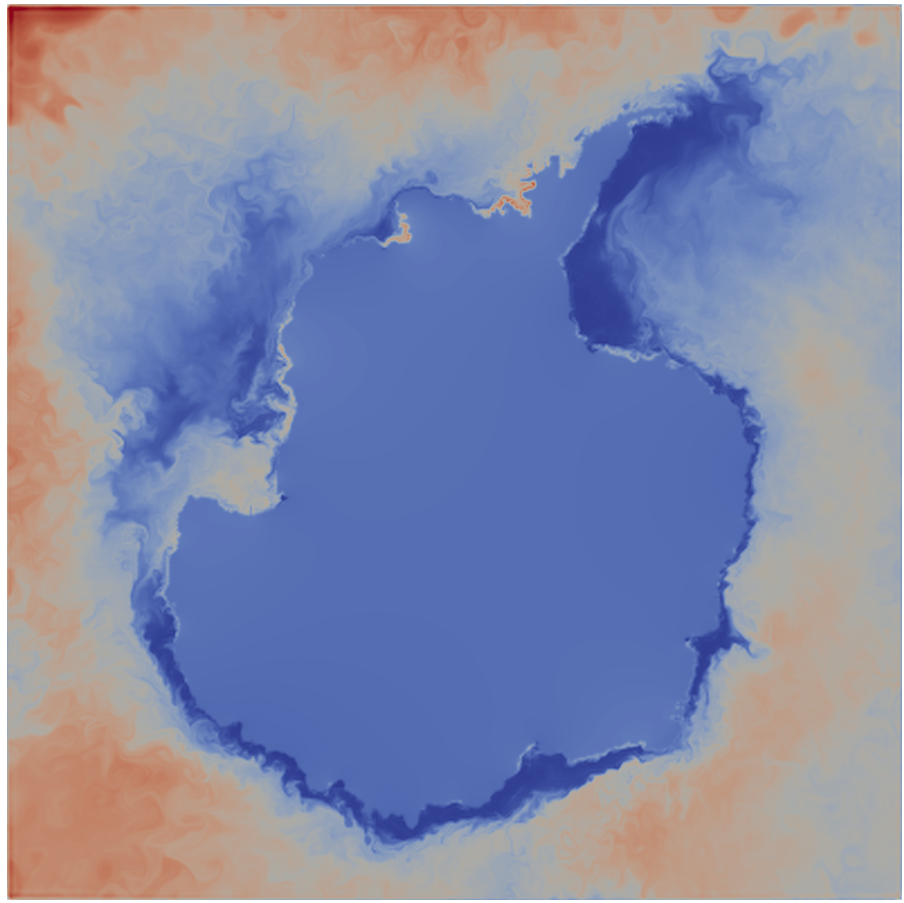}
    \includegraphics[width=\textwidth]{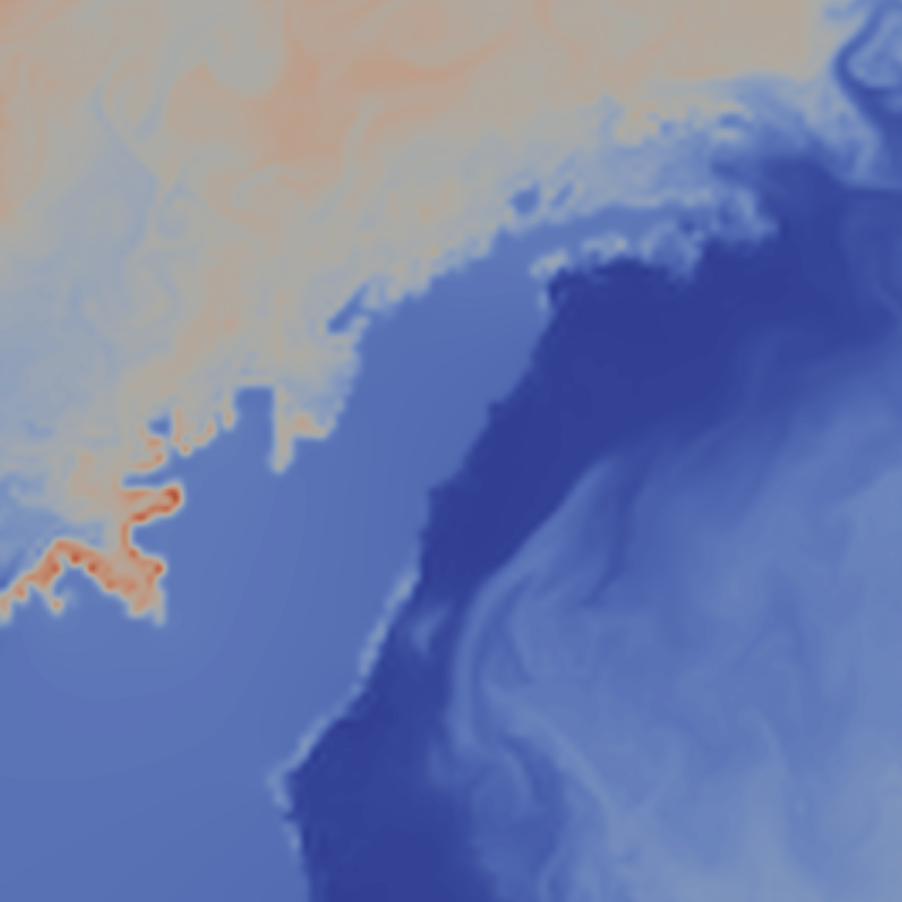}
  \end{minipage}  
  \hspace{6pt}
  \raisebox{-60pt}{\includegraphics[width=30pt]{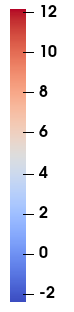}}
  \caption{Ocean temperature simulation around Antarctica. Top row: Input data on a hexagonal mesh (left), B-spline model without regularization (center), B-spline model with adaptive regularization (right). Bottom row: Detail view from the top row.} 
  \label{fig:antarctica}
\end{figure}

Figure~\ref{fig:antarctica} shows a scenario from a simulation of ocean temperatures. The data set is centered around the continent of Antarctica, for which no temperature values are given. Exhibited in the figure are two models, both of degree 2 with a 400 $\times$ 400 grid of control points. The unregularized model at center oscillates between $\pm 10^8$ along the coast (while the input data range from -2 to 10). In contrast, the regularized model (with $s^* = 5$) smoothly transitions at the coast to a near-constant value over the landmass. Regularization was performed with constraints on first and second derivatives, as described at the end of Section~\ref{sec:method}.

\begin{figure}
  \centering
  \vspace{-4ex}
  \includegraphics[width=100pt, height=100pt]{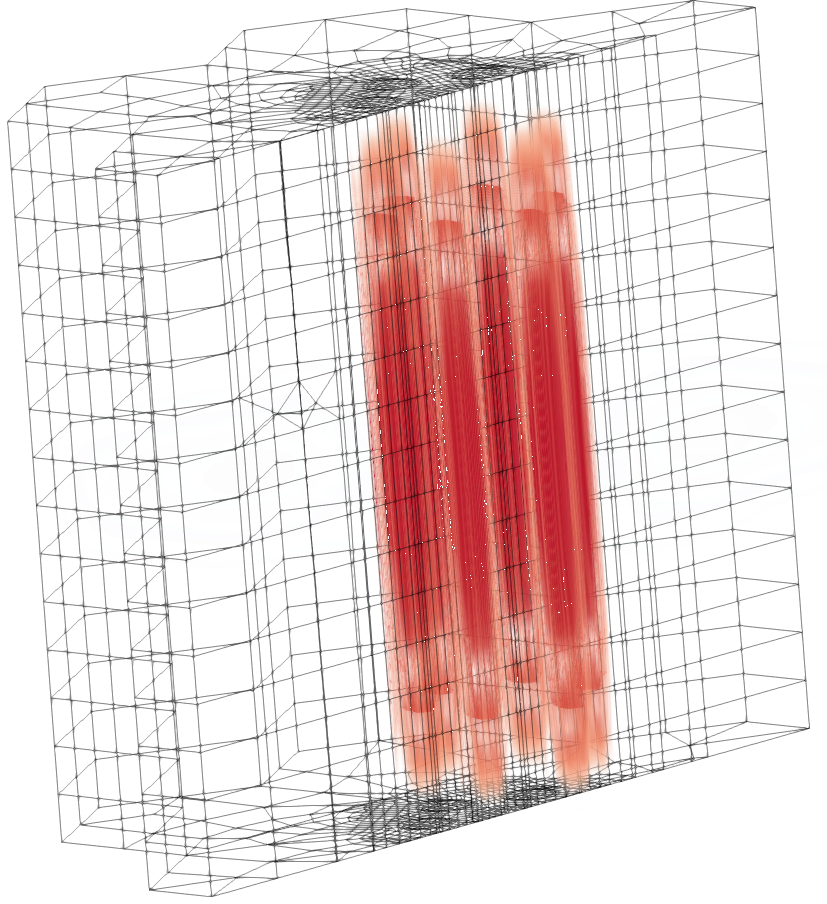}
  \includegraphics[width=100pt, height=100pt]{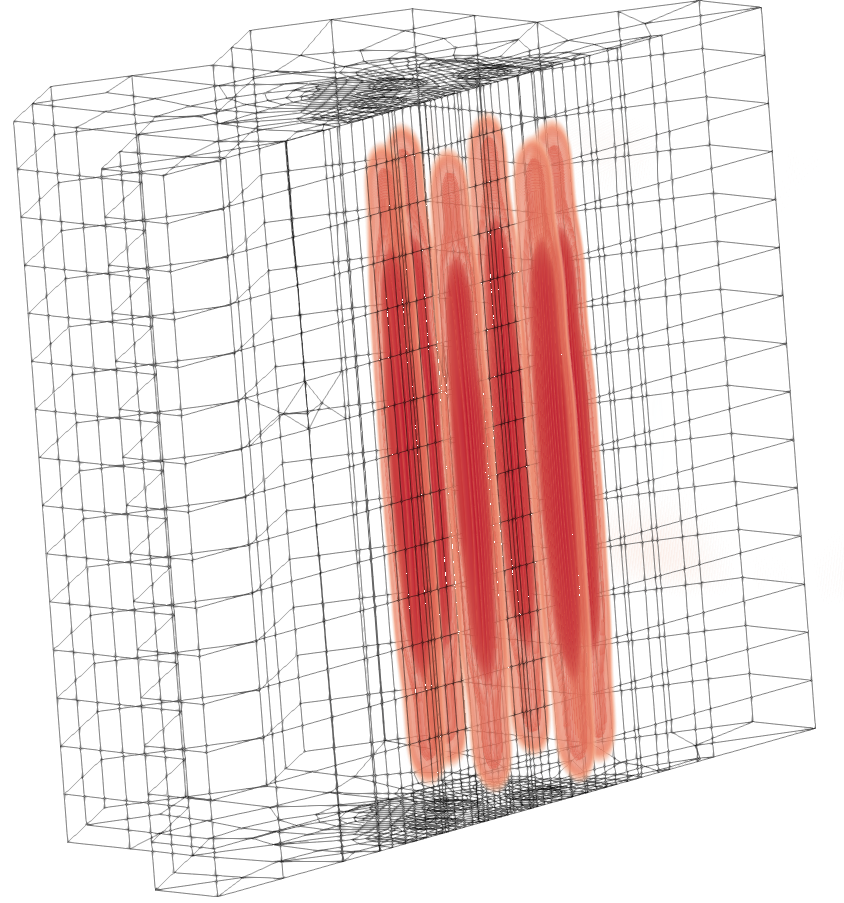}
  \includegraphics[width=100pt, height=100pt]{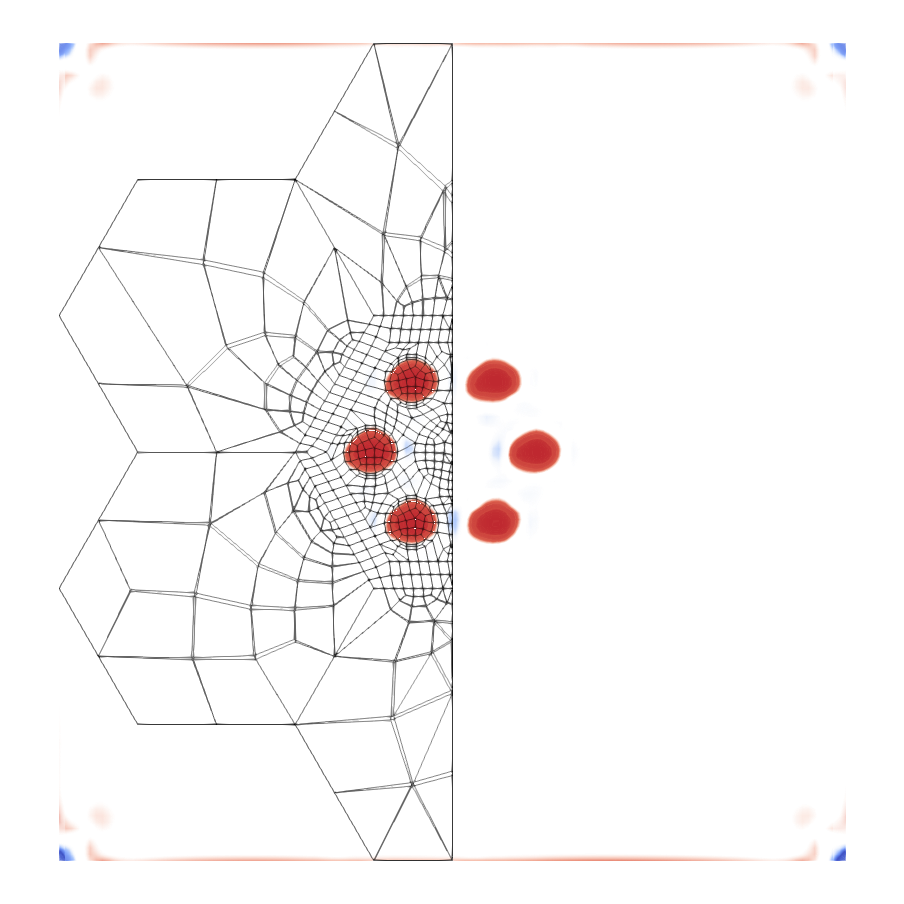}
  \raisebox{10pt}{\includegraphics[width=25pt]{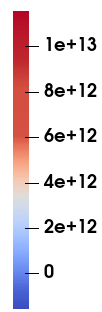}}
  \caption{Power production in a nuclear reactor simulation. Input data (left), spline model (center), spline model top view (right).}
  \label{fig:nuclear}
\end{figure}

We next consider a three-dimensional data set representing power produced in a component of a nuclear reactor (Figure~\ref{fig:nuclear}). The data are contained inside a hexagonal prism, but the B-spline model is defined on the bounding box of this prism. Hence the corners of this box are devoid of any data. 
Without regularization, the least-squares minimization does not converge, so we exhibit only our regularized model (with $s^*=10$) in Figure~\ref{fig:nuclear}. The adaptive regularization method (Figure~\ref{fig:nuclear}, center and right) produces an accurate model of the six interior ``pins'' and is well defined in the corner regions. Some artifacts are observed in the corners, but they are not significant enough to affect the interior of the model. Without regularization, the condition number of the system is infinite; with adaptive regularization, the condition number is $2.18\times10^5$. The model was produced by constraining first and second derivatives simultaneously.

\section{Future Work}\label{sec:future-work}
In the construction of our method, we imposed constraints on the first and second derivatives of the B-spline in regions where the density of input data was low or vanishing. However, these additional constraints could have been based on different orders of derivative or been unrelated to derivatives altogether. For instance, a different type of constraint would be one that penalizes deviation from a known baseline value. Another option would be to penalize B-spline values that exceed a given range (for example, the original bounds of the input data). We remark that our procedure for adaptive regularization of tensor product B-splines is separate from the type of artificial constraint imposed. Depending on the application, different constraints may be more useful, and our method allows for those to be used instead.

Another direction for future research is an investigation of the parameter $s^*$. While our method automatically varies regularization strength throughout the domain, these strengths are all relative to the parameter $s^*$. Choosing a value of $s^*$ too large can lead to overly-smoothed models. Further research into heuristics or iterative schemes to select $s^*$ automatically would allow this method to be applied with no user interaction at all.

\section{Conclusions}\label{sec:conclusions}
Modeling unstructured data sets with tensor product B-splines can be difficult due to the ill-conditioning of the fitting problem. In general, data sets with large variations in point density or regions without data exacerbate this problem to the point that artificial smoothing is necessary. However, smoothing an entire model can wash out sharp features in the data.

We introduced a regularization procedure for B-spline models that preserves features by adapting the regularization strength throughout the domain. Our method automatically varies the smoothing intensity as a function of input point density and relies on a single user-specified parameter, which we call the regularization threshold. We observe that adaptive regularization performs better than typical uniform regularization schemes that may over-smooth some regions while under-smoothing others. We also showed that our method can fit B-spline models to data sets with regions of extremely sparse point density and remain well-defined even in areas without data points. Overall, adaptive regularization of B-spline models produces smooth and accurate models for data sets which would otherwise be difficult to fit.

\section*{Acknowledgments}
This work is supported by the U.S. Department of Energy, Office of Science, Advanced Scientific Computing Research under Contract DE-AC02-06CH11357, Program Manager Margaret Lentz.

\bibliographystyle{splncs04}
\bibliography{mfa-reg}

\pagebreak
The submitted manuscript has been created by UChicago Argonne, LLC, Operator of Argonne National Laboratory (``Argonne''). Argonne, a U.S. Department of Energy Office of Science laboratory, is operated under Contract No. DE-AC02-06CH11357. The U.S. Government retains for itself, and others acting on its behalf, a paid-up nonexclusive, irrevocable worldwide license in said article to reproduce, prepare derivative works, distribute copies to the public, and perform publicly and display publicly, by or on behalf of the Government. The Department of Energy will provide public access to these results of federally sponsored research in accordance with the DOE Public Access Plan (\url{http://energy.gov/downloads/doe-public-access-plan}).
\end{document}